\documentclass[10pt,english]{amsart}
\def\margin_comment#1{\marginpar{\sffamily{\tiny #1\par}\normalfont}}
\usepackage{graphicx}
\usepackage{setspace}
\usepackage{color}
\usepackage{amsmath,amssymb,amsthm}
\usepackage{epsfig}
\usepackage{tikz}
 \usepackage{amsmath}
\usepackage{tikz-cd}
\usetikzlibrary{matrix, calc, arrows}
\usepackage{tikz}
\usetikzlibrary{arrows,chains,matrix,positioning,scopes}
\usepackage{lipsum}
\makeatletter
\tikzset{join/.code=\tikzset{after node path={%
\ifx\tikzchainprevious\pgfutil@empty\else(\tikzchainprevious)%
edge[every join]#1(\tikzchaincurrent)\fi}}}
\makeatother
\tikzset{>=stealth',every on chain/.append style={join},
         every join/.style={->}}
\tikzstyle{labeled}=[execute at begin node=$\scriptstyle,
   execute at end node=$]
\numberwithin{equation}{section} %% Comment out for sequentially-numbered
\numberwithin{figure}{section} %% Comment out for sequentially-numbered
\theoremstyle{plain}
\newtheorem*{thm*}{Theorem}
\theoremstyle{definition}
\theoremstyle{plain}
\newtheorem{thm_A}{Theorem}
\newtheorem*{defn*}{Definition}

\theoremstyle{plain}

\theoremstyle{plain} %%Delete [thm] to re-start numbering

\theoremstyle{plain}
\theoremstyle{remark}

\theoremstyle{remark}

\theoremstyle{plain}

\theoremstyle{plain}

\theoremstyle{plain}

\newtheorem*{lem*}{Lemma} %%Delete [thm] to re-start numbering
\theoremstyle{definition}

\newtheorem*{acknowledgment*}{Addentum}

\theoremstyle{plain}
\newtheorem*{ex*}{Example}
\theoremstyle{plain}

\begin{document}

\title{The Zappa-Szep product of left-orderable groups}
\author{Fabienne Chouraqui}
\maketitle
\begin{abstract}
It is  well-known that the direct product of left-orderable groups is left-orderable and that,  under a certain condition, the semi-direct product of left-orderable groups is left-orderable. We extend this result and show  that, under a  similar condition, the Zappa-Szep product of left-orderable groups is left-orderable. Moreover, we find conditions that ensure the existence of a partial left and right invariant ordering (bi-order) in  the Zappa-Szep product of bi-orderable groups and prove some properties satisfied.
\end{abstract}

%%%%%%%%%%%%%%%%%%%%%%%%%%%%%%%%%%%%%%%%%%%%%%%
\section*{Introduction}
%%%%%%%%%%%%%%%%%%%%%%%%%%%%%%%%%%%%%%%%%%%%%%%
%%%%%%%%%%%%%%%%%%%%%%%%%%%%%%%%%%%%%%%%%%%%%%%
% !TeX spellcheck = en_GB

Let $G$  be a group with subgroups $H$ and $K$ such that $G = HK$ and $H \cap K = \{1\}$. Then $G$ is isomorphic to the   Zappa-Szep product of  $H$ and $K$, denoted by $H \Join K$. If both subgroups $H$ and $K$ are normal, $G$ is their direct product $H \times K$ and if  $H$  only is normal, $G$ is their semi-direct product $H \rtimes K$. The Zappa-Szep product of groups  is a generalisation of the direct  and the semi-direct product which requires  the embedding of neither of the  factors to be normal in the product. We recall that a group $G$ is \emph{left-orderable} if there exists a  total ordering  $\prec$ of its elements which  is invariant under left multiplication, that is  $g \prec h$  implies $fg \prec fh$ for all $f,g,h$ in $G$. It is  well-known that the direct product of left-orderable groups is left-orderable and that,  under a certain condition, the semi-direct product of left-orderable groups is left-orderable \cite{kopytov}, \cite{book-order}[p.27]. We extend this result and show  that, under a  similar condition, the Zappa-Szep product of left-orderable groups is left-orderable. Before we state our first result, we introduce some notations.
\begin{defn*}
Let $H$, $K$ be groups. Let $\alpha$  be the homomorphism defined by $\alpha:K \rightarrow \operatorname{Sym}(H) $, $k \mapsto \alpha_k$, where $k \in K$ and  $\operatorname{Sym}(H)$ denotes the group of bijections of $H$. Let  $\beta:H \rightarrow \operatorname{Sym}(K)$, $h \mapsto\beta_h$, be such that  $\beta_{h_1h_2}=\beta_{h_2} \circ \beta_{h_1}$ and $\beta_{1}=Id_{K}$ ($\beta$ an anti-homomorphism).
Assume $(\alpha,\beta)$ satisfies: $\alpha_{k}(h_1h_2)=\,\alpha_{k}(h_1) \,\alpha_{\beta_{h_1}(k)}(h_2)$    and  
$\beta_{h}(k_1k_2)=\,\beta_{\alpha_{k_2}(h)}(k_1)\,\beta_{h}(k_2)$. \emph{The Zappa-Szep product $H \Join K$ of $H$ and  $K$ with respect to the pair $(\alpha,\beta)$} is  the set $H \times K$ endowed with the following product:
\[(h_1,k_1)\,(h_2,k_2)=\, (h_1\,\alpha_{k_1}(h_2),\,\beta_{h_2}(k_1)\,k_2)\] 
The identity is $(1,1)$ and the inverse of an element $(h,k)$ is
\[(h,k)^{-1} = (\alpha_{k^{-1}}(h^{-1}),\,\beta_{h^{-1}}(k^{-1})\,)\]
In fact, $H\times \{1\}$ and $K\times \{1\}$  are subgroups of $H \Join K$,  isomorphic to $H$ and $K$, respectively.
\end{defn*}
The  Zappa-Szep product is also called crossed product,  bi-crossed product,  knit product, two-sided semidirect product. We refer to \cite{brin}, \cite{gigel}, \cite{szep}, \cite{zappa} for details. Note that every element in  $H \Join K$ is uniquely represented by $hk$, with $h \in H$ and $k \in K$,  and $kh$ is equal to  $\alpha_{k}(h)\,\beta_{h}(k)$ in $H \Join K$. 
\begin{thm_A}\label{thm_intro_1}
Let $G$ be the  Zappa-Szep product of the groups $H$  and $K$ with respect to $(\alpha,\beta)$. Assume $H$  and $K$  are left-orderable groups.  Let $P_H$ and $P_K$ denote the   positive cones (of a left order) of $H$  and $K$ respectively. Assume the following condition $(*)$ is satisfied: \[\alpha_{k}(P_H) \subseteq P_H, \; \forall k \in K\]
Then, there exists a  total left order  $\prec$ on  $G= H \Join K$, with positive cone $P$, such that (the embedding of ) $K$ is a convex subgroup with respect to  $\prec$ and $P_K=P\cap K$. 
 \end{thm_A}
We introduce some definitions and refer to \cite{remtulla}, \cite{deroin}, \cite{levi}, \cite{book-order}, \cite{kopytov}, \cite{linnel},  \cite{rolfsen}, \cite{sikora}. Let $G$ be a left-orderable group with  a  strict total left order  $\prec$. An element $g$, $g \in G$,  is called \emph{positive} if $1 \prec g$ and the set of positive elements $P$ is called the \emph{positive cone of} $\prec$. The positive cone $P$ satisfies: \\
$(1)$ $P$ is a semigroup, that is $P\cdot P \subseteq P$\\
$(2)$ $G$ is partitioned, that is $G= P \cup P^{-1} \cup \{1\}$ and  $P \cap P^{-1}=\emptyset$\\
 Conversely, if there exists a subset $P$ of $G$ that satisfies $(1)$ and $(2)$, then $P$ determines a unique total  left order $\prec$ defined by $g \prec h$ if and only if $g^{-1}h \in P$. If $P$ satisfies only $(1)$, the left order obtained is partial. A subgroup $N$ of a left-orderable group $G$ is called \emph{convex} if for any $x,y,z \in G$ such that $x,z \in N$ and $x \prec y \prec z$, we have $y \in N$. 
\begin{proof}\emph{of Theorem $1$}
Let $P \subseteq G$ be defined by:  $g=hk \, \in P$ if $h\in P_H$ or if $h=1$, $k \in P_K$. We show there exists  a total left order $\prec$  on $G= H \Join K$ with positive cone $P$. First, we  prove  $P$ is a semigroup. Let $g=hk \, \in P$ and $g'=h'k' \, \in P$. If $h'=1$ and $h=1$, then $k,k' \in P_K$ and $kk'\, \in P_K$, since $P_K$ is a semigroup, so $gg' \in P$. If  $h'=1$ and $h \neq 1$, then $h \in P_H$ and  so $gg'=hkk' \in P$. If $h' \neq 1$, then $h' \in P_H$ and  $gg'=\,hkh'k' =\, h\,\alpha_k(h')\beta_{h'}(k)\,k'$. From the assumption $(*)$, $\alpha_k(h') \in P_H$,  so if $h=1$, then $gg'=\, \alpha_k(h')\beta_{h'}(k)\,k' \,\in P$ and if $h \neq 1$,  $ h\,\alpha_k(h')\in P_H $ also,  since $P_H$ is a semigroup, so $gg' \in P$. Next, we prove that  given $g=hk \neq 1$ in $G= H \Join K$, either $g$ belongs to $P$ or $g^{-1}=\alpha_{k^{-1}}(h^{-1})\,\beta_{h^{-1}}(k^{-1})$ belongs to $P$. Assume  $h=1$. 
If  $k \in P_K$, then $g \in P$, otherwise $k^{-1}\in P_K$, since $P_K$ partitions $K$, and  then   
$g^{-1}=k^{-1} \in P$. Furthermore, $g$ cannot belong to $P \cap P^{-1}$, since it would contradict $P_K \cap P_K^{-1}=\emptyset$. Assume  $h\neq 1$. If  $h \in P_H$, then $g \in P$, 
otherwise $h^{-1}\in P_H$, since $P_H$ partitions $H$, and  from  $(*)$,  $\alpha_{k^{-1}}(h^{-1}) \in P_H$. So,  
$g^{-1}=\alpha_{k^{-1}}(h^{-1})\,\beta_{h^{-1}}(k^{-1})$ belongs to $P$. Assume  $g \in P \cap P^{-1}$, then $h \in P_H$ and  $\alpha_{k^{-1}}(h^{-1}) \in P_H$. It holds that $\alpha_{k^{-1}}(h^{-1}) = (\alpha_{\beta_{h^{-1}}(k^{-1})}(h))^{-1}$. Indeed, 
on one side $\alpha_{k^{-1}}(1)=1$ and on the other side $\alpha_{k^{-1}}(1)=\,\alpha_{k^{-1}}(h^{-1}h)=\,\alpha_{k^{-1}}(h^{-1}) \,\alpha_{\beta_{h^{-1}}(k^{-1})}(h)$.   From  $(*)$,  $\alpha_{\beta_{h^{-1}}(k^{-1})}(h)$ belongs to $P_H$, so $\alpha_{k^{-1}}(h^{-1}) = (\alpha_{\beta_{h^{-1}}(k^{-1})}(h))^{-1}$ belongs to $P_H \cap P_H^{-1}$ and this is a contradiction. So, $P$ satisfies the conditions 
$(1)$ and $(2)$ and it determines uniquely a total left order $\prec$.\\
Let $g=hk \in G$ and $k' \in K$.  Assume $1 \prec g \prec k'$ and assume by contradiction that $h \neq 1$. From  $1 \prec g$,  we have $h \in P_H$. From $ hk \prec k'$, we have $k'^{-1}hk \prec 1$. But, 
$k'^{-1}hk =\, \alpha_{k'^{-1}}(h)\,\beta_{h}(k'^{-1})\,k\,\succ 1$, since, from $(*)$, $\alpha_{k'^{-1}}(h)$ belongs to $P_H$. So, $h=1$.
\end{proof}
 In the following example, we  remark it is sometimes useful to consider  a semi-direct product  as a special case of  Zappa-Szep product  in order to show it is left-orderable. Indeed, this permits to check condition $(*)$ either on a positive cone of  $H$ or on a positive cone of $K$. Let $H= \operatorname{Gp}\langle y,z \mid y^3=z^3 \rangle$.  Since $H= \mathbb{Z}\ast_{3\mathbb{Z}}\mathbb{Z}$ and $\mathbb{Z}$ is (left) orderable,  $H$ is left-orderable \cite{book-order}[p.178]. Let $K= \operatorname{Gp}\langle a,b \mid aba=bab \rangle$, the braid group on $3$ strands, $K$ is left-orderable \cite{dehornoy}. Let  $\alpha:K \rightarrow \operatorname{Sym}(H)$  be the trivial homomorphism and let  $\beta:H \rightarrow \operatorname{Sym}(K)$,  the anti-homomorphism defined by $\beta_y=\beta_z=(a,b)$.  With respect to $(\alpha,\beta)$, let  $G= \,H \Join K$   in this specific order; $G$ is presented by  $\operatorname{Gp}\langle y,z ,a,b \mid aba=bab ,\,y^3=z^3\,,ay=yb, by=ya,az=zb,bz=za \rangle$.  We show $G$ is  left-orderable, using Thm.\ref{thm_intro_1}. Each element $h \in H$ admits a normal formal  $h=\,y^{\epsilon_1}z^{\mu_1}y^{\epsilon_2}z^{\mu_2}...y^{\epsilon_m}z^{\mu_m}\,\Delta^{n_h}$, where $\Delta=y^3=z^3$ is central in $H$,  $-1\leq \epsilon_i,\mu_i \leq 1$ and $n_h \in \mathbb{Z}$. We define  $exp(h)=\epsilon_1 +\mu_1+..+\epsilon_m+\mu_m+3n_h$ and  $P_H$ to be the set of elements $h \in H \setminus\{1\}$ such that $exp(h) >0$, or if $exp(h)=0$, $n_h>0$, or if $exp(h)=n_h=0$, $\epsilon_1 \neq 0$. The set $P_H$ is a positive cone (the proof appears in the appendix) and  $P_H$ satisfies trivially  the condition  $(*)$ from Thm.\ref{thm_intro_1}. So, $G$ is  left-orderable. 
 
 Note that a group is left-orderable if and only if it is right-orderable. This is well illustrated here with the existence of a 
 symmetric version of this result. Indeed,  if  $H$  and $K$  are right-orderable groups,  with $Q_H$ and $Q_K$ positive cones of right orders of $H$  and $K$  satisfying $\beta_{h}(Q_K) \subseteq Q_K, \; \forall h \in H$. Then there exists a  total right order  $<$ on  $G= H \Join K$, with positive cone $Q$ (defined by $g=hk \, \in Q$ if $k\in Q_K$ or if $k=1$, $h \in Q_H$), such that (the embedding of ) $H$ is a convex subgroup with respect to  $<$.  We say a  group $G$ is \emph{partially (totally) bi-orderable} if there exists a  partial (total)  ordering  $\ll$ of its elements which  is invariant under left and right multiplication, that is  $g \ll h$  implies $fgk \ll fhk$ for all $f,g,h,k$ in $G$. In particular, a set $P$ determines a partial  bi-order if and only if  $P$  is a semigroup, and satisfies $gPg^{-1}\subseteq P$ for all $g \in G$;  $P$ determines a total  bi-order if additionally  $G= P \cup P^{-1} \cup \{1\}$ and  $P \cap P^{-1}=\emptyset$. A natural question is when the  Zappa-Szep product of  bi-orderable groups is a bi-orderable group. In the following theorem, we give conditions that ensure the existence of a partial bi-order.

\begin{thm_A}\label{thm_intro_2}
Let $G$ be the  Zappa-Szep product of the groups $H$  and $K$ with respect to $(\alpha,\beta)$. Assume $H$  and $K$  are partially  bi-orderable groups.  Let $P_H$ and $P_K$ denote the   positive cones (of a partial bi-order) of $H$  and $K$ respectively. Assume the following conditions are satisfied: \\$(*)$  $\alpha_{k}(P_H) \subseteq P_H, \; \forall k \in P_K,$ and 
$\beta_{h}(P_K) \subseteq P_K, \; \forall h \in P_H$ 
 \\$(**)$ $k\,P_H\,k^{-1} \subseteq P_H, \; \forall k \in K,$ and $h\,P_K\,h^{-1} \subseteq P_K, \; \forall h \in H$\\
Then, there exists a  partial bi-order  $\ll$ on  $G= H \Join K$. Furthermore, if $h,h' \in H$ satisfy $h\ll  h'$ then $\alpha_k(h) \,\ll \alpha_k(h')$ and $\beta_h(k) \,\ll \beta_{h'}(k)$,  $\forall k \in K$  and if $k,k' \in K$ satisfy $k\ll  k'$ then $\alpha_k(h) \,\ll \alpha_{k'}(h)$ and $\beta_h(k) \,\ll \beta_{h}(k')$,  $\forall h \in H$.  
 \end{thm_A}
\begin{proof}
Let $P \subseteq G$ be defined by:  $g=hk \, \in P$ if $h\in P_H$ and  $k \in P_K$; if $h=1$,  $k \in P_K$ or if $k=1$, $h\in P_H$. We show $P$ is a semigroup and  $gPg^{-1}\subseteq P$ for all $g \in G$. Let $g=hk \, \in P$ and $g'=h'k' \, \in P$, then  $gg'=\,hkh'k' =\, h\,\alpha_k(h')\beta_{h'}(k)\,k'$. From the assumption $(*)$, $\alpha_k(h') \in P_H$ and $\beta_{h'}(k) \in P_K$,  so  $h\,\alpha_k(h')\in P_H $ and $\beta_{h'}(k)\,k' \in P_K$, since $P_H$ and  $P_K$ are semigroups, so $gg' \in P$. Next, let  $g=hk \, \in G$ and $h'k' \in P$, we show $gh'k'g^{-1} \in P$: $gh'k'g^{-1} = hkh'k'k^{-1}h^{-1}=\,(h\,(kh'k^{-1})\,h^{-1})\,(h\,(kk'k^{-1})\,h^{-1})$. From $(**)$, $h'\in P_H$  implies $kh'k^{-1} \in P_H$ and so $h\,(kh'k^{-1})\,h^{-1} \in P_H$, since $P_H$ is the positive cone of a bi-order. The element $h\,(kk'k^{-1})\,h^{-1} \in P_K$, first from  $kk'k^{-1} \in P_K$, and next from $(**)$.\\
 Let  $h,h' \in H$ satisfy $h\ll  h'$, and let $k \in K$, then $kh\ll  kh'$, that is  $\alpha_k(h)\beta_h(k) \,\ll\alpha_k(h')\beta_{h'}(k)$. 
 Since $\ll$ is a bi-order, this implies $1\ll(\alpha_k(h))^{-1}\,\alpha_k(h')\,\beta_{h'}(k)(\beta_h(k) )^{-1}$ and from the definition of $\ll$, we have $(\alpha_k(h))^{-1}\,\alpha_k(h')\,\in P_H$ and $\beta_{h'}(k)(\beta_h(k) )^{-1} \in P_K$, that is  
  $\alpha_k(h) \,\ll \alpha_k(h')$ and $\beta_h(k) \,\ll \beta_{h'}(k)$. In case $h=1\, h'\neq 1$,  this means $\alpha_{k}(P_H) \subseteq P_H, \; \forall k \in K$ and  $h' \gg 1$ implies $\beta_{h'}(k)\gg k$, $\forall k \in K$; in case $h\neq 1\, h'=1$,  this means that $h \ll 1$ implies $\alpha_k(h)\ll 1$ and  $\beta_{h}(k)\ll k$, $\forall k \in K$. The same proof works for the symmetric statement. 
\end{proof}
 It would be interesting to know if there exists a natural extension of total bi-orders of $H$ and $K$ that defines  a total bi-order on $H \Join K$ and $H \Join K$ is not a direct nor a semi-direct product.
\section*{Appendix}
 We recall each element $h \in H$ admits a normal formal  $h=\,y^{\epsilon_1}z^{\mu_1}y^{\epsilon_2}z^{\mu_2}...y^{\epsilon_m}z^{\mu_m}\,\Delta^{n_h}$, where $\Delta=y^3=z^3$ is a central element in $H$,  $-1\leq \epsilon_i,\mu_i \leq 1$ and $n_h \in \mathbb{Z}$. We define  $exp(h)=\epsilon_1 +\mu_1+..+\epsilon_m+\mu_m+3n_h$ and  $P_H$ to be the set of elements $h \in H \setminus\{1\}$ such that $exp(h) >0$ (class (A)), or if $exp(h)=0$, $n_h>0$ (class (B)), or if $exp(h)=n_h=0$, $\epsilon_1 \neq 0$ (class (C)). We show $P_H$ is a positive cone. First, we show $P_H$ is a semigroup.
  Let $h=\,y^{\epsilon_1}z^{\mu_1}...y^{\epsilon_m}z^{\mu_m}\,\Delta^{n_h}$, $h'=\,y^{\epsilon'_1}z^{\mu'_1}...y^{\epsilon'_m}z^{\mu'_m}\,\Delta^{n_{h'}}$, with  $-1\leq \epsilon_i,\mu_i, \epsilon'_i,\mu'_i \leq 1$, be elements of $P_H$. We need to check several cases in order to show $hh' \in P_H$. \\
  Case $1$: if $h$ is in class (A) and $h'$ is in class (A) or (B) or (C), then  $hh'$ is also in (A).\\
  Case $2$: if $h, h'$ are in class (B), then  $hh' \in P_H$, since $exp(hh')=0$ and $n_{hh'} \geq n_h +n_{h'} -1>0$. Indeed, there are only two cases in which $n_{hh'}$ decreases: if $h$ ends with $y^{-1}$ and $h'$ begins with $y^{-1}$, or if  $h$ ends with $z^{-1}$ and $h'$ begins with $z^{-1}$, and in both cases $n_{hh'} =n_h +n_{h'} -1$.
  Case $3$: assume  $h, h'$ are in class (C). Note that if $h=\,y^{\epsilon_1}z^{\mu_1}...y^{\epsilon_m}z^{\mu_m}\,\Delta^{n_h}$, $-1\leq \epsilon_i,\mu_i \leq 1$,  satisfies $exp(h)=n_h=0$, then   the length of $h$ is necessarily even. Furthermore, if $h$ is  in class (C), then $h=\,y^{\pm 1}z^{\mu_1}y^{\epsilon_2}z^{\mu_2}...y^{\epsilon_m}z^{\pm 1}$ and the same holds for $h'$. So, $hh'$ is also in (C).\\
  Case $4$: if $h$ is in class (B) and $h'$ is in class  (C), then $h=\,y^{\epsilon_1}z^{\mu_1}...y^{\epsilon_m}z^{\mu_m}\,\Delta^{n_h}$, $n_h>0$
and $h'=\,y^{\pm 1}z^{\mu'_1}...y^{\epsilon'_m}z^{\pm 1}$. It holds that $n_h -1 \leq n_{hh'} \leq n_h +1$ and we need only to check the case $n_{hh'}=n_h -1$. If $n_h -1 > 0$, then $hh'$ is in class (B). If $n_h -1 = 0$, then $h=\,y^{\epsilon_1}z^{\mu_1}...z^{\mu_{m-1}}y^{-1}\,\Delta$, 
 $h'=\,y^{- 1}z^{\mu'_1}...y^{\epsilon'_m}z^{\pm 1}$ and $hh'= \,y^{\epsilon_1}z^{\mu_1}...z^{\mu_{m-1}}\,y\,z^{\mu'_1}...y^{\epsilon'_m}z^{\pm 1}$. We show $hh'$ is in class  (C), that is $\epsilon_1 \neq 0$. On one-hand, $exp(h)=\epsilon_1+\mu_1+...+\mu_{m-1}\,-1+3$ and on the second-hand $exp(h)=0$. By induction on $m$, it holds that $\epsilon_1+\mu_1+...+\epsilon_{m-1}+\mu_{m-1}=-2$, with $-1\leq \epsilon_i,\mu_i \leq 1$, implies $\epsilon_1 \neq 0$. \\ 
 Case $5$: if $h$ is in class (C) and $h'$ is in class  (B), then $h=\,y^{\pm _1}z^{\mu_1}...y^{\epsilon_m}z^{\pm 1}$, 
   and $h'=\,y^{\epsilon'_1}z^{\mu'_1}...y^{\epsilon'_m}z^{\mu'_m}\,\Delta^{n_{h'}}$, $n_{h'}>0$. It holds that
    $n_{h'} -1 \leq n_{hh'} \leq n_{h'} +1$ and in case $n_{hh'}=n_{h'}-1=0$,  $hh'$ is necessarily in class  (C), as $hh'$ begins with $y^{\pm _1}$.\\
      Next, we show $P_H$ partitions $H$, that is given $h \in H \setminus\{1\}$, either $h \in P_H$ or $h^{-1} \in P_H$.
  If $exp(h)>0$, then  $h \in P_H$, otherwise assume $exp(h) \leq 0 $. If  $exp(h)<0$, then  $h^{-1} \in P_H$, since  $exp(h^{-1})=\,-exp(h)$.
  So, assume  $exp(h)=0$. If $n_h>0$, then  $h \in P_H$, otherwise assume $n_h \leq 0 $. If $n_h < 0$, then $h^{-1} \in P_H$, since  $n_{h^{-1}}=\,-n_{h}$. So, assume $n_h = 0$. If $h=\,y^{\epsilon_1}z^{\mu_1}...y^{\epsilon_m}z^{\mu_m}\,\Delta^{n_h}$, $-1\leq \epsilon_i,\mu_i \leq 1$,  satisfies $exp(h)=n_h=0$, then   there are two possibilities: either $h=\,y^{\pm 1}z^{\mu_1}y^{\epsilon_2}z^{\mu_2}...y^{\epsilon_m}z^{\pm 1}$  or $h=\,z^{\pm 1}y^{\epsilon_2}z^{\mu_2}...y^{\epsilon_{m-1}}z^{\mu_{m-1}}y^{\pm 1}$. If $h=\,y^{\pm 1}z^{\mu_1}y^{\epsilon_2}z^{\mu_2}...y^{\epsilon_m}z^{\pm 1}$, then $h \in P_H$. Otherwise, if  $h=\,z^{\pm 1}y^{\epsilon_2}z^{\mu_2}...y^{\epsilon_{m-1}}z^{\mu_{m-1}}y^{\pm 1}$, then $h^{-1}=y^{\mp 1}z^{-\mu_{m-1}}y^{-\epsilon_{m-1}}z^{-\mu_2}y^{-\epsilon_2}...z^{\mp 1}$, that is $h^{-1} \in P_H$.

%%%%%%%%%%%%%%%%%%%%%%%%%%%%%%%%%%%%%%%%%%%%%%%%%%%%%%%%%%%%%%%%
\bigskip\bigskip\noindent
{ Fabienne Chouraqui,}

\smallskip\noindent
University of Haifa at Oranim, Israel.
\smallskip\noindent\\
E-mail: {\tt fabienne.chouraqui@gmail.com,  fchoura@sci.haifa.ac.il}


\begin{thebibliography}{40}
\bibitem{szep} J. Szep, \emph{On the structure of groups which can be represented as the product of two subgroups},
Acta Sci. Math. Szeged {\bf 12} (1950), 57-61.
\end{thebibliography}
\end{document}